\documentclass[a4paper,12pt,thmsa]{amsart}
\usepackage[a4paper,marginratio={1:1},scale={0.72,0.74},footskip=7mm,headsep=10mm]{geometry}

\usepackage[utf8]{inputenc}
\usepackage[english]{babel}

\usepackage{amsfonts}
\usepackage{amssymb,amsmath,latexsym}
\usepackage{graphicx,subfigure}
\usepackage[dvips]{color}
\usepackage[T1]{fontenc}
\usepackage[active]{srcltx}
\usepackage{amsmath}
\usepackage{amsfonts}
\usepackage{amssymb}
\usepackage{hyperref}
\usepackage{psfrag}
\usepackage{color}
\usepackage{url}
\usepackage{amsthm}
\usepackage{array}
\usepackage{pst-tree}
\usepackage{lscape}
\usepackage{dsfont}
\usepackage[T1]{fontenc}
\usepackage{pstricks,pstricks-add}\usepackage{mathrsfs}  
\usepackage{wrapfig}
\usepackage[labelformat=empty]{caption}
\title[A lifetime of excursions through random walks and L\'evy processes]
{A lifetime of excursions through random walks and L\'evy processes}

\author{Lo\"ic Chaumont}

\address{Lo\"ic Chaumont, Univ Angers, CNRS, LAREMA, SFR MATHSTIC, F-49000 Angers, France}

\email{loic.chaumont@univ-angers.fr}

\author{Andreas Kyprianou}

\address{Andreas Kyprianou -- 
Department of Mathematical Sciences, University of Bath, Claverton Down, BA2 7AY}

\email{a.kyprianou@bath.ac.uk}

\begin{document}

\begin{abstract} 
In celebration of Professor Ron Doney’s 80th birthday, we provide a summary of his academic career and contributions to probability theory, as one of the UK’s leading probabilists for over 50 years. A version of this note also serves as an introductory article to the volume in honor of Ron Doney [Birkh\"auser, 2021], which includes an additional 17 research papers produced by Ron's colleagues and friends. 
\end{abstract}

\maketitle


Through this volume, it is with the greatest of admiration that we  pay tribute to  the mathematical achievements of Professor Ron Doney. His career has spanned generations of
probabilists and his work continues to play a significant role in the community. In addition to the  major contributions he has made   in the theory of random walks and L\'evy processes, Ron is equally appreciated for the support he has given to younger colleagues.
The sentiment and desire to organise both a workshop and this volume to honour his lifetime  achievements surfaced naturally a couple of years ago at the L\'evy process meeting in Samos, as it became apparent that Ron was approaching his 80th birthday. A huge appreciation of his standing in the community meant it was very easy to find willing participation in both projects. As a prelude to this {\it hommage \`a Ron}, let us spend a little time reflecting on his career and his main achievements. 

\bigskip

Ron grew up in working-class Salford in the North West of  Greater Manchester at a time when few leaving school would attend university. 
Having a  veracious appetite for reading, Ron  spent long hours  as a schoolboy in Manchester Central Library where he cultivated his intellect. Coupled with an obvious   talent for mathematics, he found his way to the University of Durham. There he studied Mathematics as an undergraduate and  continued all the way to a PhD under the supervision of Harry Reuter, graduating in 1964 aged just 24. Ron's thesis entitled {\it `Some problems on random walks'} was no doubt inspired by the shift  in interests of Harry Reuter from analysis to probability. 
At the time of his PhD,  all doctoral activity in probability theory in the UK was essentially supervised either by Harry Reuter or David Kendall. Because of this, Reuter and Kendall  formed StAG, the {\it Stochastic Analysis Group}, which would meet regularly under the auspices of the London Mathematical Society. 
 From a very young age, Ron thus had the opportunity to engage with his contemporaries. Back in the 1960s they were  precious few in number compared with vast numbers of the probability PhD community present today in the UK and  included the likes of    Rollo Davidson, David Williams, Daryl Daley, David Vere-Jones, Nick Bingham, John Hawkes and John Kingman.

\bigskip

The time that Ron completed his PhD coincided with a period of  expansion in the UK higher education system which proved to offer numerous opportunities during the first years of his academic career. From his PhD, Ron was successfully appointed directly to a lectureship at the University of East Anglia. He spent the academic year 1964-65 there, but quickly moved on to lectureship at Imperial College London in 1965, coinciding  with the appointment of Harry Reuter  to a chair there. During this early phase, Ron had a slow start to his publication record. His first two papers \cite{1966-1, 1966-2} concern random walks in three dimensions, followed by a paper concerning higher dimensional version of the renewal theory \cite{1966-3}. 

\bigskip

By 1970,  Ron had moved back to his native  Manchester. He  joined the then world-famous Manchester-Sheffield School of Probability and Statistics, formally taking up a lectureship at Manchester's Statistical Laboratory. At this point in time he moved away from his initial work on random walks to the theory of Galton-Watson processes and, what were then called, general branching processes (today they would rather come under the heading of Crump-Mode-Jagers (CMJ) processes). Although seemingly a change in direction, this was a very natural move for anyone who harboured interests in random walks and renewal processes. Indeed, whilst the theory of branching processes has become significantly more exotic in recent years, the interplay of these two fields still remains highly pertinent today. Concurrently with the work of Peter Jagers, Ron produced a cluster of articles through to the mid 1970s looking at growth properties of CMJ processes in which, among other things, he demonstrated the central role that renewal processes play; \cite{1971, 1972-1, spine, 1973, 1974, 1975, 1976, 1977-1, 1977-2}. As early as 1972, one finds computations in his work which echo what would later become known as the method of {\it spines}; \cite{spine}. Here, Ron also made contributions to underlying functional equations and the so-called $x\log x$ condition, which precede a number of similar results in  the setting of more general spatial models such as  branching random walks. The papers, \cite{1974, 1975} with Nick Bingham, are also interesting to reflect upon in terms of how collaborations of the day were conducted. In a period of no internet or email, and with the probability community in the UK being very few in number and  thinly spread, Ron maintained communication with Nick through hand-written and mailed letters during the epoch of their overlapping interests. In a process that has largely been replaced by Google, this involved sharing sample calculations, summaries of articles they had found and  broader mathematical ideas. In the case of Ron and Nick, this led not only two these two papers, but calculations that lay dormant, surfacing over a decade  later in \cite{1988}. 

\begin{wrapfigure}[18]{r}{0.4\textwidth}
  \begin{center}
\includegraphics[height=6cm]{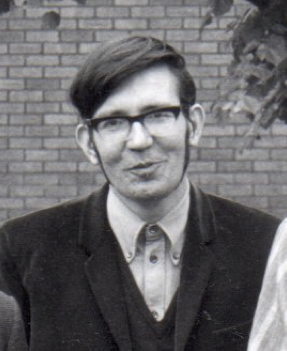}
  \end{center}
\caption{\it Ron Doney in 1972, Manchester.}
\end{wrapfigure}

\bigskip

By 1977, Ron was back to random walks, albeit now in one dimension. As a passive observer, it appears that there was an awakening in Ron's understanding of how many  problems still remained open for general one-dimensional random walks, particularly when looking at them in terms of excursions. 
Aside from some works on Markov chains and Brownian motion, \cite{1984-2, 1984-1,  1989-4},  Ron's  work focused mostly with random walks  through the 1980s. The prominence of his contributions lay with his use of Spitzer's condition as well as the problem of characterising the moments and tail behaviour of ladder variables in relation to assumptions on the underlying random walk; \cite{1980-1, 1980-2, 1982-1, 1982-2, 1983, 1989-1, 1989-2, 1989-3}. He began a growing interest in arcsine laws and random walks in the domain of attraction of stable laws as well as stable processes themselves; \cite{1985, 1987, WHF, 1988}. Most notable of the latter is a remarkable paper on the Wiener-Hopf factorisation of stable processes \cite{WHF}, the significance of which would become apparent many  years later after 2010 thanks to continued work of a number of authors, most notably  Alexey Kuznetsov. It was also during that 1980s that Ron took two sabbaticals to Canada.  The first stay was in Vancouver in 1980-81, visiting Cindy Greenwood, which also allowed him the chance to connect with Sidney Port, Ed Perkins and John Walsh. For the  second, he was visiting George O'Brien in York University, Toronto in 1988-89.

\bigskip

Moving into the 1990s, many things changed for Ron's research, least of all, his rate of publication, \cite{1991-1, 1991-2, 1991-levy, 1992-1, 1992-2, 1993-1, 1993-2, 1993-3, J1, J2, J3, J4, 1995, J5, 1996, 1997, J7, 1998-1, 1998-2, 1998-3, 1998-4, Al}. The previous work he had done on ladder heights and Spitzer's condition culminated in one of Ron's most important and widely appreciated results: For random walks (and shortly after for L\'evy processes), he proved  that Spitzer's condition was equivalent to the convergence of the positivity probability, \cite{1995, 1997}. Ron found himself catapulted into a rapidly growing and much better organised global community of probabilists with shared interests to his own. Although he largely continued publishing in the context of random walks, it was during this decade that Ron became increasingly exposed to the theory of L\'evy processes. 

\begin{wrapfigure}[15]{r}{0.4\textwidth}
  \begin{center}
\includegraphics[scale=1]{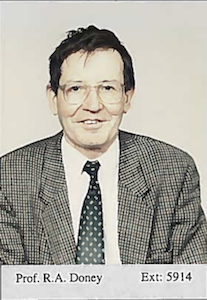}
  \end{center}
\caption{\it Ron Doney’s departmental photo as it appeared on the entrance to the Mathematics Tower in the 1990s.}
\end{wrapfigure}

\bigskip

The 1990s  was also the decade that saw Ron begin to jointly co-author many more of his papers. There were a number of factors at play here. During the first twenty five years of Ron's career it was quite common for academics to author papers alone, least of all because, within the field of probability theory, probabilists were few and far between. The 1990s was the decade of globalisation and the prevalence of email and internet made international connectivity and communication much easier.  But perhaps most importantly in this respect was Ron's collaboration with a young French mathematician by the name of Jean Bertoin, which dramatically opened up his relationship with the Parisian (and more generally the French) school 
of probability for the rest of his career. 

\bigskip

In the early 1990s, Ron had made contact with Jean Bertoin, who was very interested  in random walks, but  also in capturing and expanding on the less-well explored  theory of  L\'evy processes. It was Ron's first paper on L\'evy processes, \cite{1991-levy}, which had caught the interest of Jean and they began sending each other preprints. They first met in person  in Luminy during one of the {\it Journ\'ees de Probabilit\'es} in 1992 organised by Az\'ema and Yor. 
At the time, Ron was interested in a paper by Keener on the simple RW conditioned to stay positive that had been published that year. Although others had written on the topic of conditioned random walks, it was the joint work of Jean and Ron \cite{J4} that ensued, which formalised  a robust approach to the notion of such conditionings and led the way for a number of articles on conditioned processes, particularly in the setting of L\'evy processes. Reflecting on earlier remarks about how the theory of random walks and branching processes are so closely intertwined, it is worth noting that the paper \cite{J4}  ended up playing a hugely influential role in the theory of  branching random walks, branching Brownian motion and (growth) fragmentation processes, where conditioning of the {\it spine} in the spirit of their work proved to be instrumental in understanding the so-called {\it derivative martingale}. 

\bigskip

The collaboration between Jean and Ron was relatively intensive for a period of 3-4 years with six of their seven co-authored papers, \cite{J1, J2, J3, J4, J5, J6, J7}, appearing between 1994 and 1997. During this period, Jean would visit Ron  regularly, often working at his home in the village of Whaley Bridge, South East of Manchester. Their daily working routine  would be interlaced with regular  hikes out into  the  Peak District which lay just beyond Ron's back door.

\bigskip

Through his collaboration with Jean Bertoin, Ron also started to visit Paris more often. There,  he became acquainted with a new generation of young probabilists who were being guided towards the theory of L\'evy processes. Lo\"{\i}c Chaumont was the first PhD student of Jean who also became a long-time collaborator of Ron. Their collaboration spans 8 papers to date, \cite{Ch1, Ch2, Ch3, Ch4, Ch5, Ch6, Ch7, Ch8}, in which they cover the study of perturbed Brownian motion, distributional decompositions of the general Wiener-Hopf factorisation and L\'evy processes conditioned to stay positive. It is perhaps the latter, \cite{Ch5, Ch6}, for which they are best known as co-authors, building on the PhD thesis of  Lo\" {\i}c for conditioned L\'evy processes that had, in turn, grown out of the  formalisation for conditioned random walks that Ron had undertaken with Jean Bertoin. It was also during the mid 1990s that Larbi Alili, a contemporary of Lo\"{\i}c and PhD student of Marc Yor,  became the postdoc of Ron on a competitively funded EPSRC  project; cf \cite{Ch1, Al, Alili}. Philippe Marchal, another gifted young probabilist from the Paris school, was also a regular  visitor to Manchester during this period. 
Another young French probabilist whose work greatly impressed Ron at the  turn of the Millennium  was Vincent Vigon. 
Aside from Ron's admiration of Vincent's unexpected emergence from  Rouen rather than Paris,  what impressed him the most was that Vigon had established a necessary and sufficient condition, in the form of an explicit  integral test, for when a  L\'evy process of unbounded variation with no Gaussian component creeps; in particular, this result showed that Ron's previous conjecture on this matter, which had been assimilated from the long-term behaviour of random walks,  was wrong. 

\bigskip

The new Millennium brought about yet further change for Ron.  His publication rate went up yet another gear, with almost as many articles published during this decade as in the previous two. This was all  thanks to his increased exposure to collaborative partnership as well as the inevitable depth of understanding of random walks and L\'evy processes he had acquired; \cite{Ch2, Ch3, 2001, Alili, 2001-3, 2001-4, R1, R2, 2003, 2003-2, 2004, R3, 2004-2, R4, 2004-4, Ch4, Tusheng, 2005, 2005-2, R5, Ch5, R6, 2006, 2006-2, R7, 2007-2, R8, J7, Ch6, 2008, R10, 2009-2}. It was also during this decade that Ron began an extremely fruitful collaboration with Ross Maller, publishing 10 papers together; \cite{R1, R2, R3, R4, R5, R6, R7, R8, J7, R10, R11}.  On a visit to the UK, Ross was advised by Charles Goldie to go and spend time in Manchester visiting Ron. He did and they immediately started producing material.
In a number of important papers \cite{R3, R2, R10}, Ross and Ron first investigated the 
asymptotic behaviour of random walks and L\'evy processes at deterministic times and at  first passage times across a fixed level. In another series of remarkable works \cite{R6, R7, R8, J7}  they considered first 
passage times across  power law boundaries of random walks, L\'evy processes and their reflected version 
at the infimum.  In particular they obtained necessary and sufficient conditions for these first passage times to 
have finite moments. 

\bigskip

By now, there was widespread renewed interest in the theory of fluctuations for 
L\'evy processes and the study of their overshoots had  become very popular. Two articles by Ron, written with Phil Griffin \cite{2003}, \cite{2004-4} bear witness to this. It was also during this time that Ron collaborated with 
Andreas Kyprianou and  wrote one of his most cited articles on the so-called {\it quintuple law}. The latter gives a distributional identity for a suite of five important and  commonly used path functionals  of a L\'evy process at first passage over a fixed level  in terms of its ladder potentials and L\'evy  measure. In essence, the result  constitutes a `disintegration' of the Wiener--Hopf factorisation. 
Ron also began a fruitful collaboration with  
Mladen Savov during this period; \cite{R10, 2010, 2010-2}. Mladen, who came to Manchester from Buglaria, proved to be the most accomplished of Ron's several PhD students.

\bigskip

A very important event of this decade for Ron was his invitation to give a lecture at the famous Saint-Flour summer 
school in 2005. It was arguably the first major recognition of Ron's career by the mathematical community. 
For Ron, this  also presented the opportunity to write a book on fluctuation theory for L\'evy processes, \cite{2007-2}, which remains an important reference 
in the domain to date. 
 Ron found himself centre-stage as part of a huge community of researchers now working specifically in the field of random walks and L\'evy processes. He spoke at many venues, including a rapid succession of workshops and congresses devoted to L\'evy processes and was elected as a Fellow of the Institute of Mathematical Statistics in 2006. Moreover, in 2005, Manchester hosted the 4th international workshop on L\'evy processes. This was a huge undertaking given the large number of attendees, but nonetheless an important moment that asserted the importance of Manchester's probability group and, in particular, Ron's identity as a highly accomplished researcher in this field.

\bigskip

\begin{wrapfigure}[18]{r}{0.5\textwidth}
\begin{center}
\includegraphics[height=5cm]{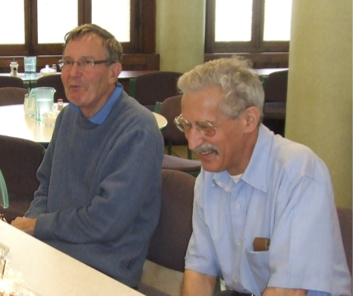}
\end{center}
\caption{\it Taken on the 24th July 2009 at the 8th International Conference on L\'evy processes, held in Angers, Ron Doney sits with Marc Yor to celebrate Marc's 60th birthday.}
\end{wrapfigure}

The next decade, 2010-2020, saw an invitation  by the Bernoulli Society for Ron to give a prestigious
plenary lecture at the 2013 Stochastic Processes and their Applications congress in Boulder, Colorado.
Together with his lectures at  Saint Flour, this stands as  quite an important recognition of his achievements from the international mathematical community. Technically speaking, Ron retired in 2006, however, as a special case given the quality of his research output, Manchester had  gladly  continued his appointment beyond his 67th birthday,  right through until 2014, on a purely research basis. 
Now in his 70s, one might say Ron worked at a slower rate, however, although publishing in lower volume,  one sees remarkable crowning quality in his work with an exceptional number of his papers appearing in  the top two probability journals {\it Annals of Probability} and {\it Probability Theory and Related Fields};
\cite{2010, Ch7, 2010-2, 2011, 2012, 2012-2, 2012-3, 
2012-4, 2013, 2013-2, 2015, R11, 2016, 2018, 2019, 2020, Ch8}. 

\bigskip

These works are marked by two  projects 
in which Ron revisited old obsessions. The first is a cluster of papers, co-authored with V\'{\i}ctor Rivero. Ron went back to the first passage times 
of L\'evy processes and obtained sharp results for the local behaviour of their distribution, \cite{2013}, \cite{2015}, 
\cite{2016}. The second concerns   improving the classical Gnedenko and Stone local limit theorems. In collaboration with Francesco Caravenna \cite{2019}, Ron obtained necessary and sufficient conditions for random walks in the domain of attraction of a 
stable law to satisfy the strong renewal theorem. This work, as well as \cite{2020}, solve a long-standing problem
which dates back to the 1960's. This second achievement is the one that  Ron himself quietly admits  he is most proud of, and rightly so.

\bigskip

As Ron's career  wound down, he had the pleasure of watching the probability group in Manchester dramatically grow in size, with 10 members around the time of his final retirement, something he had dreamed of for many years. Ron was appointed in Manchester at a time that it stood as a global  stronghold for  probability and statistics. The Manchester Statistical Laboratory was  half of the  very unique two-university Manchester-Sheffield School of Probability and Statistics,  the brain child of Joe Gani. 
With changing times the Manchester-Sheffield School disbanded and,  aside from  Fredos Papangelou, who joined in 1973, Ron was the only other probabilist who remained in the Statistical Laboratory for a number of years to come. The 1980s were hard times for British academia, but the growth in the international community around Ron's research  interests through the 1990s was mirrored in the growth of the probability group in Manchester. The turn of the  Millennium saw the merger of the University of Manchester (or rather, more formally, the Victoria University of Manchester) with UMIST (University of Manchester Institute of Science and Technology), which opened the door to  new opportunities.
Many of the probability appointments in Manchester since then clearly reflect the strong  association of  Manchester with the theory of random walks and L\'evy processes; something that is  directly tethered to Ron's towering achievements as a researcher. 

\bigskip

As alluded to at the start of this article, Ron is appreciated as much for his encouragement of young researchers as he is for the mathematics that he has produced. Just as the authors of this article see Ron as one of the major influencing characters in their careers, both through academic mentorship and mathematical  discourse,  so do many others among our community, both in the UK and around the globe. There are simply too many to list here, moreover,  an attempt to do so  would  carry the risk that we forget names. But it should be said that, when the idea of holding a workshop for Ron's 80th birthday surfaced, this was carried forward by an emotional surge of support from the many who belong to the aforementioned list. 

\bigskip

As a humble researcher who cares little for the limelight, Ron did not  always get the honours he deserved.  
Towards his retirement, the number of researchers in probability theory exploded exponentially and Ron's contribution to a classical field, which   now lies  in the DNA of many modern research endeavours, is often overlooked.
One goal of   this volume and the accompanying conference is to try to correct this.

\bigskip

We hope the contents of this volume will stimulate Ron to think about his next piece of work. He currently has 99 publications and we are all keenly awaiting his 100th!

\end{document}